\documentclass[12pt, reqno]{amsart}
\allowdisplaybreaks
\begin{document}
\setcounter{page}{1}

\title[\hfilneg \hfil time scales in three variables]
{ Some new dynamic Inequality  on time scales in three variables}

\author[Deepak B. Pachpatte\hfil \hfilneg]
{Deepak B. Pachpatte}

\address{Deepak B. Pachpatte \newline
 Department of Mathematics,
 Dr. Babasaheb Ambedkar Marathwada University, Aurangabad,
 Maharashtra 431004, India}
\email{pachpatte@gmail.com}


\begin{abstract}
 In this paper we obtain the estimates on some dynamic integral inequalities in three variables which can be used to study certain dynamic equations. We give some applications to convey the importance of our result.
\end{abstract}

\maketitle

\section{Introduction}
   The study of time scales was initiated in 1989 by Stefan Hilger \cite{HIG}  in his Ph.D dissertation. Since then many authors have studied the dynamic inequalities on time scales.  Some analytic inequalities on time scales in one and two variables is studied in \cite{Hus, Ozk, Sun, Yeh} by various authors. The authors in  \cite{Boh3,Dbp1,Dbp2,Dbp3} have obtained some interesting dynamic integral and iterated inequalities on time scales.
  Motivated by the results  above in this  paper we establish new explicit bounds on some dynamic inequalities in three variables  which are useful in solving certain dynamic equations.

In what follows $\mathbb{R}$ denotes the set of real numbers, I=[a,b] and $\mathbb{T}$ denotes arbitrary time scales. We say that $f:\mathbb{T} \to \mathbb{R}$ is rd-continuous provided $f$ is continuous right dense point of $\mathbb{T}$ and has a  finite left sided limit at each left dense point of  $\mathbb{T}$ and will be denoted by $C_{rd}$. Let $\mathbb{T}_1$ and $\mathbb{T}_2$ be two time scales with atleast two points and $\Omega  = \mathbb{T}_1  \times \mathbb{T}_2$ and $H = \Omega  \times I$.
The basic information about time scales can be found in \cite{Boh1,Boh2}.
 Now we give the Lemma given in \cite{Fer} which is required in proving our result.

Lemma [\cite{Fer}] Let $u,a,f \in C'_{rd} \left( {\mathbb{T}_1  \times \mathbb{T}_2 ,\mathbb{R}_ +  } \right)$ and $a$ is nondecreasing in each of the variables. If
\[
u\left( {x,y} \right) \le a\left( {x,y} \right) + \int\limits_{x_0 }^x {\int\limits_{y_0 }^y {f\left( {s,t} \right)u} } \left( {s,t} \right)\Delta t\Delta s,
\tag{1.1}\]
for $\left( {x,y} \right) \in \mathbb{T}_1  \times \mathbb{T}_2 $ then
\[
u\left( {x,y} \right) \le a\left( {x,y} \right)e_{\int\limits_{y_0 }^y {f\left( {x,t} \right)\Delta t} } \left( {x,x_0 } \right),
\tag{1.2}\]
for $\left( {x,y} \right) \in \mathbb{T}_1  \times \mathbb{T}_2 $.

\section{Main Results}
Now we give our main result in the following theorem

 \textbf{Theorem 2.1} Let $u,p,q,f \in C_{rd} \left( {H,\mathbb{R}_ +  } \right)$ and $c \ge 0$ be a constant. If
 \[
u\left( {x,y,z} \right) \le p_1 \left( {x,y,z} \right) + p_2 \left( {x,y,z} \right)\int\limits_{x_0 }^x {\int\limits_{y_0 }^y {\int\limits_a^b {f\left( {s,\tau ,q} \right)} } } u\left( {s,\tau ,q} \right)\Delta q\Delta \tau \Delta s,
\tag{2.1}\]
for $(x,y,z) \in H$, then
\[
u\left( {x,y,z} \right) \le p_1 \left( {x,y,z} \right) + p_2 \left( {x,y,z} \right)C\left( {x,y} \right)e_{Q\left( {x,y,z} \right)} \left( {x,x_0 } \right),
\tag{2.2}\]
where
\[
Q\left( {x,y,z} \right) = \int\limits_{y_0 }^y {\int\limits_a^b {f\left( {s,\tau ,q} \right)p_2 \left( {s,\tau ,q} \right)\Delta q\Delta \tau \Delta s} },
\tag{2.3}\]
\[
C\left( {x,y} \right) = \int\limits_{x_0 }^x {\int\limits_{y_0 }^y {\int\limits_a^b {f\left( {s,\tau ,q} \right)p_1 \left( {s,\tau ,q} \right)\Delta q\Delta \tau \Delta s} } }.
\tag{2.4}\]
\textbf{Proof.} Now let
   \[
M\left( {s,t} \right) = \int\limits_a^b {f\left( {s,\tau ,q} \right)p_2 \left( {s,\tau ,q} \right)\Delta q}.
\tag{2.5}\]
Then $(2.1)$ becomes
\[
u\left( {x,y,z} \right) \le p_1 \left( {x,y,z} \right) + p_2 \left( {x,y,z} \right)\int\limits_{x_0 }^x {\int\limits_{y_0 }^y {M\left( {s,\tau } \right)} } \Delta \tau \Delta s.
\tag{2.6}\]
Now put
\[
W\left( {x,y} \right) = \int\limits_{x_0 }^x {\int\limits_{y_0 }^y {M\left( {s,\tau } \right)} } \Delta \tau \Delta s.
\tag{2.7}\]
Then $W(x,y_0)=W(x_0,y)=0$
and
\[
u\left( {x,y,z} \right) \le p_1 \left( {x,y,z} \right) + p_2 \left( {x,y,z} \right)W\left( {x,y} \right).
\tag{2.8}\]
From $(2.7),(2.5),(2.8)$ we have
 \begin{align*}
&W^{\Delta _1 \Delta _2 } \left( {x,y} \right) \\
&= M\left( {x,y} \right) \\
& = \int\limits_a^b {f\left( {x,y,q} \right)u\left( {x,y,q} \right)\Delta q}  \\
& \le \int\limits_a^b {f\left( {x,y,q} \right)\left[ {p_1 \left( {x,y,z} \right) + p_2 \left( {x,y,z} \right)W\left( {x,y} \right)} \right]\Delta q}  \\
&= W\left( {x,y} \right)\int\limits_a^b {f\left( {x,y,q} \right)u\left( {x,y,q} \right)\Delta q}  + \int\limits_a^b {f\left( {x,y,q} \right)p_1 \left( {x,y,q} \right)\Delta q}  \\
& = \int\limits_a^b {f\left( {x,y,q} \right)p_2 \left( {x,y,q} \right)\Delta q}  + \int\limits_a^b {f\left( {x,y,q} \right)p_1 \left( {x,y,q} \right)\Delta q}.
\tag{2.9}
 \end{align*}
Now from $(2.9)$ above we have by taking delta integral
\begin{align*}
 W^{\Delta _1 } \left( {x,y} \right)
 &\le \int\limits_{y_0 }^y {\int\limits_a^b {W\left( {x,\tau } \right)f\left( {x,\tau ,q} \right)} } p_2 \left( {x,\tau ,q} \right)\Delta q\Delta \tau  \\
 & + \int\limits_{y_0 }^y {\int\limits_a^b {f\left( {x,\tau ,q} \right)} } p_1 \left( {x,\tau ,q} \right)\Delta q\Delta \tau.
\tag{2.10}
 \end{align*}
Again delta integrating above $(2.10)$ we have
 \begin{align*}
 W\left( {x,y} \right)
 &\le \int\limits_{x_0 }^x {\int\limits_{y_0 }^y {\int\limits_a^b {W\left( {s,\tau } \right)f\left( {x,\tau ,q} \right)} } p_2 \left( {x,\tau ,q} \right)\Delta q\Delta \tau }  \\
 &+ \int\limits_{x_0 }^x {\int\limits_{y_0 }^y {\int\limits_a^b {f\left( {s,\tau ,q} \right)} } p_1 \left( {s,\tau ,q} \right)\Delta q\Delta \tau }.
\tag{2.11}
 \end{align*}
Put

$B\left( {x,y} \right) = \int\limits_a^b {f\left( {x,y,q} \right)p_2 \left( {x,y,q} \right)\Delta q} $,

and

$C\left( {x,y} \right) = \int\limits_{x_0 }^x {\int\limits_{y_0 }^y {\int\limits_a^b {f\left( {s,\tau ,q} \right)} } p_1 \left( {s,\tau ,q} \right)\Delta q\Delta \tau }$.

We get from $(2.11)$
\[
W\left( {x,y} \right) \le \int\limits_{x_0 }^x {\int\limits_{y_0 }^y {B\left( {s,\tau } \right)W\left( {s,\tau } \right)} \Delta \tau \Delta s}  + C\left( {x,y} \right).
\tag{2.12}\]

Clearly $C(x,y)$ is nondecreasing in $\Omega$ then applying Lemma to $(2.12)$,  we get
\[
W\left( {x,y} \right) \le C\left( {x,y} \right)e_{\overline Q \left( {x,y} \right)} \left( {x,x_0 } \right),
\tag{2.13}\]
where
\[
\overline Q \left( {x,y} \right) = \int\limits_{y_0 }^y {B\left( {x,\tau } \right)} \Delta \tau.
\tag{2.14}\]
Now using $(2.13)$ in
\[
u\left( {x,y,z} \right) \le p_1 \left( {x,y,z} \right) + p_2 \left( {x,y,z} \right)W\left( {x,y} \right),
\]
we get the result $(2.2)$.

\section{Applications}
Now in this section we give some applications of our results.
 Consider the dynamic integral equation of the form
\[
u\left( {x,y,z} \right) = g\left( {h,y,z} \right) + \int\limits_{x_0 }^x {\int\limits_{y_0 }^y {\int\limits_a^b {F\left( {x,y,z,s,t,q} \right)} } } \Delta q\Delta t\Delta s,
\tag{3.1}\]
for $(x,y,z) \in H$  where $g \in C_{rd} \left( {H,\mathbb{R}} \right)$, $F \in C_{rd} \left( {H^2  \times \mathbb{R},\mathbb{R}} \right)$.

Now our next theorem deals with the estimate of solution of $(3.1)$.

 \textbf{Theorem 3.1} Suppose the function $F$ in $(3.1)$ satisfy the conditions
 \[
\left| {F\left( {x,y,z,s,t,q,u} \right)} \right| \le r\left( {x,y,z} \right)f\left( {s,t,q} \right)\left| u \right|,
\tag{3.2}\]
where $r,f \in C_{rd} \left( {H,\mathbb{R}} \right)$. If $u(x,y,z)$ is a solution of equation $(3.1)$ then
\[
\left| {u\left( {x,y,z} \right)} \right| \le \left| {g\left( {h,y,z} \right)} \right| + r\left( {x,y,z} \right)C_2 \left( {x,y,z} \right)e_{Q_2 \left( {x,y,z} \right)} \left( {x,x_0 } \right),
\tag{3.3}\]
where
\[
C_2 \left( {x,y,z} \right) = \int\limits_{x_0 }^x {\int\limits_{y_0 }^y {\int\limits_a^b {f\left( {s,t,q} \right)} } } \left| {g\left( {s,t,q} \right)} \right|\Delta q\Delta t\Delta s,
\tag{3.4}\]
\[
Q_2 \left( {x,y,z} \right) = \int\limits_{y_0 }^y {\int\limits_a^b {f\left( {x,t,q} \right)} } r\left( {x,t,q} \right)\Delta q\Delta t,
\tag{3.5}\]
for $(x,y,z) \in H$.

\textbf{Proof.}  Let $u \in C_{rd} \left( {H,\mathbb{R}} \right)$ be a solution of $(3.1)$, we have
\begin{align*}
\left| {u\left( {x,y,z} \right)} \right|
&\le \left| {g\left( {h,y,z} \right)} \right| + \int\limits_{x_0 }^x {\int\limits_{y_0 }^y {\int\limits_a^b {\left| {F\left( {x,y,z,s,t,q} \right)} \right|} } } \Delta q\Delta t\Delta s \\
&\le \left| {g\left( {h,y,z} \right)} \right| \\
&+ r\left( {x,y,z} \right)\int\limits_{x_0 }^x {\int\limits_{y_0 }^y {\int\limits_a^b {f\left( {s,t,q} \right)} } } \left| {u\left( {s,t,q} \right)} \right|\Delta q\Delta t\Delta s.
\tag{3.6}
\end{align*}
Now applying the Theorem $(2.1)$ gives the estimate $(3.3)$.

Now for obtaining estimates in our next theorem we suppose that $F$ satisfies Lipschitz type conditions.

 \textbf{Theorem 3.2} Suppose that the function $F$ in $(3.1)$ satisfies the condition
\[
\left| {F\left( {x,y,z,s,t,q,u} \right) - F\left( {x,y,z,s,t,q,v} \right)} \right| \le r\left( {x,y,z} \right)f\left( {s,t,q} \right)\left| {u - v} \right|,
\tag{3.7}\]
where $r,f \in C_{rd} \left( {H,\mathbb{R}} \right)$. If $u(x,y,z)$ is a solution of $(3.1)$ then
\begin{align*}
&\left| {u\left( {x,y,z} \right) - g\left( {x,y,z} \right)} \right| \\
&\le k\left( {x,y,z} \right) + r(x,y,z)C_3 (x,y,z)e_{Q_2 \left( {x,y,z} \right)} \left( {x,x_0 } \right),
\tag{3.8}
\end{align*}
for $(x,y,z) \in H$ where
\[
C_3 \left( {x,y,z} \right) = \int\limits_{x_0 }^x {\int\limits_{y_0 }^y {\int\limits_a^b {f\left( {s,t,q} \right)\left| {k\left( {s,t,q} \right)} \right|\Delta q\Delta t\Delta s} } },
\tag{3.9}\]

and
\[
k\left( {x,y,z} \right) = \int\limits_{x_0 }^x {\int\limits_{y_0 }^y {\int\limits_a^b {\left| {F\left( {x,y,z,s,t,q,g(s,t,q)} \right)} \right|} } } \Delta q\Delta t\Delta s,
\tag{3.10}\]
for $(x,y,z) \in H$.

\textbf{Proof.} Let $u \in C_{rd} \left( {H,\mathbb{R}} \right)$ be a solution of equation $(3.1)$. Then we have

\begin{align*}
&\left| {u\left( {x,y,z} \right) - g\left( {x,y,z} \right)} \right| \\
&\le \int\limits_{x_0 }^x {\int\limits_{y_0 }^y {\int\limits_a^b {\left| {F\left( {x,y,z,s,t,q,u(s,t,q)} \right)} \right|} } } \Delta q\Delta t\Delta s \\
& \le \int\limits_{x_0 }^x {\int\limits_{y_0 }^y {\int\limits_a^b {\left| {F\left( {x,y,z,s,t,q,u(s,t,q)} \right)} \right.} } }  \\
&\left. { - F\left( {x,y,z,s,t,q,u(s,t,q)} \right)} \right|\Delta q\Delta t\Delta s \\
& + \int\limits_{x_0 }^x {\int\limits_{y_0 }^y {\int\limits_a^b {\left| {F\left( {x,y,z,s,t,q,g(s,t,q)} \right)} \right|} } } \Delta q\Delta t\Delta s \\
& \le k\left( {x,y,z} \right) \\
& + r\left( {x,y,z} \right)\int\limits_{x_0 }^x {\int\limits_{y_0 }^y {\int\limits_a^b {f\left( {s,t,q} \right)} } } \left| {u\left( {s,t,q} \right) - h\left( {s,t,q} \right)} \right|\Delta q\Delta t\Delta s,
\tag{3.11}
\end{align*}
for $(x,y,z) \in H$.

Now an application of theorem $(2.1)$ to $(3.11)$ gives the estimate $(3.8)$.

Now we consider equation $(3.1)$ and the integral equation
\[
h\left( {x,y,z} \right) = v\left( {x,y,z} \right) + \int\limits_{x_0 }^x {\int\limits_{y_0 }^y {\int\limits_a^b {G\left( {x,y,z,s,t,q,h(x,y,z)} \right)\Delta q\Delta t\Delta s} } },
\tag{3.12}\]
for $v \in C_{rd} \left( {H,\mathbb{R}} \right)$,$G \in C_{rd} \left( {H^2  \times \mathbb{R},\mathbb{R}} \right)$.

Now we give the following theorem.

 \textbf{Theorem 3.3}
Suppose the function $F$ in $(3.1)$ satisfies the condition $(3.7)$ then for every solution $h \in C_{rd} \left( {H,\mathbb{R}} \right)$ of $(3.11)$ and $u \in C_{rd} \left( {H,\mathbb{R}} \right)$ solution of $(3.1)$ we have the estimates
\begin{align*}
\left| {u\left( {x,y,z} \right) - h\left( {x,y,z} \right)} \right|
&\le \left[ {\overline g \left( {x,y,z} \right) + \overline k \left( {x,y,z} \right)} \right] \\
& + r\left( {x,y,z} \right)C_4 \left( {x,y,z} \right)e_{Q_2 \left( {x,y,z} \right)} \left( {x,x_0 } \right).
\tag{3.13}
\end{align*}
for $(x,y,z) \in H$ in which
\[
C_4 \left( {x,y,z} \right) = \int\limits_{x_0 }^x {\int\limits_{y_0 }^y {\int\limits_a^b {f\left( {s,t,q} \right)\left[ {\overline g \left( {s,t,q} \right) + \overline k \left( {s,t,q} \right)} \right]\Delta q\Delta t\Delta s} } }.
\tag{3.14}\]
\[
\overline g \left( {x,y,z} \right) = \left| {g\left( {x,y,z} \right) - v\left( {x,y,z} \right)} \right|.
\tag{3.15}\]
\begin{align*}
\overline k \left( {x,y,z} \right)
&= \int\limits_{x_0 }^x {\int\limits_{y_0 }^y {\int\limits_a^b {\left| {F\left( {x,y,z,s,t,q,h(s,t,q)} \right)} \right.} } }  \\
&\left. { - G\left( {x,y,z,s,t,q,h(s,t,q)} \right)} \right|\Delta q\Delta t\Delta s.
\tag{3.16}
\end{align*}
for $(x,y,z) \in H$.

 \textbf{Proof.} Since $u(x,y,z)$ and $v(x,y,z)$ are respectively solutions of $(3.1)$ and $(3.12)$ we have
\begin{align*}
&\left| {u\left( {x,y,z} \right) - h\left( {x,y,z} \right)} \right| \\
&\le \left[ {\overline g \left( {x,y,z} \right) + \overline k \left( {x,y,z} \right)} \right] \\
&+ \int\limits_{x_0 }^x {\int\limits_{y_0 }^y {\int\limits_a^b {\left| {F\left( {x,y,z,s,t,q,u(s,t,q)} \right)} \right.} } }  \\
&\left. { - F\left( {x,y,z,s,t,q,h(s,t,q)} \right)} \right|\Delta q\Delta t\Delta s \\
&+ \int\limits_{x_0 }^x {\int\limits_{y_0 }^y {\int\limits_a^b {\left| {F\left( {x,y,z,s,t,q,u(s,t,q)} \right)} \right.} } }  \\
&\left. { - G\left( {x,y,z,s,t,q,h(s,t,q)} \right)} \right|\Delta q\Delta t\Delta s \\
&\le \overline g \left( {x,y,z} \right) + \overline k \left( {x,y,z} \right) \\
&+ r\left( {x,y,z} \right)\int\limits_{x_0 }^x {\int\limits_{y_0 }^y {\int\limits_a^b {f\left( {s,t,q} \right)\left| {u\left( {s,t,q} \right) - h\left( {s,t,q} \right)} \right|\Delta q\Delta t\Delta s} } }.
\tag{3.17}
\end{align*}
Now an application of Theorem $2.1$ to $(3.17)$ yields $(3.13)$.

\end{document}